\documentclass[11pt]{article}

 \usepackage{amsfonts,amssymb,graphicx}
\usepackage[colorlinks=true,linkcolor=blue,citecolor=blue]{hyperref}

\newtheorem{lemm}{Lemma}
\newtheorem{theo}{Theorem}

\newcommand{\C}[1][]{\ensuremath{{\mathbb{C}^{#1}} }}
\newcommand{\R}[1][]{\ensuremath{{\mathbb{R}^{#1}} }}

\renewcommand{\S}[1][]{\ensuremath{{\mathbb{S}^{#1}} }}

\newcommand{\<}{\langle}
\renewcommand{\>}{\rangle}
\newcommand{\ga}{\gamma}

\newcommand{\al}{\alpha}

\newcommand{\be}{\beta}

\date{}
\def\mail#1#2{{\tt #1}@{\tt #2}}
\author{Henri Anciaux\footnote{The author is supported by SFI (Research Frontiers Program)}}
\title{Two non existence results for the self-similar equation in Euclidean 3-space}
\begin{document}
\maketitle
 \vspace{-2.4em} \centerline{ \small email:
\mail{henri.anciaux}{staff.ittralee.ie}}

\vspace{ 3em}

\centerline{\textbf {Abstract}}

\smallskip

{\small We prove that the only  self-similar surfaces of Euclidean
$3$-space which are foliated by circles are  the self-similar
surfaces of revolution discovered by S. Angenent and that the only
ruled, self-similar surfaces are the cylinders over planar
self-similar curves.}

\medskip

\centerline{\small \em 2000 MSC: 53C44\em }







\section*{Introduction}
The Mean Curvature
Flow (denoted by MCF in the following) is the gradient flow of the area functional on the space of
$n$-submanifolds of some Riemannian manifold. From the viewpoint
of analysis, this flow is governed by a non-linear parabolic
equation. Although classical results of analysis show short-time
existence of the MCF, understanding its long-time behaviour is a
hard problem which requires to control
 the possible singularities that may appear
along the flow.

Self-similar flows arise as special solutions of the MCF
 that preserve the shape of the evolving
submanifold. Analytically speaking, this amounts to making a
particular Ansatz in the parabolic PDE describing the flow in
order to eliminate the time variable and reduce the equation to an
elliptic one.

 The simplest and most important example of a self-similar flow is when the
evolution is a homothety. Such a self-similar submanifold $X$ with
mean curvature vector $\vec{H}$ satisfies the following
non-linear, elliptic system:
 $$  \vec{H} + \lambda X^{\perp} =0,  $$
where $X^{\perp}$ stands for the projection of the position vector
$X$ onto the normal space. If $\lambda$ is any strictly positive
constant, the submanifold shrinks in finite time to a single point
under the action of the MCF, its shape remaining
 unchanged. If $\lambda$ is strictly negative, the submanifold will expand,
its shape again remaining the same; in this case the submanifold
is necessarily non-compact. The case of vanishing $\lambda$ is the
well-known case of a minimal submanifold, which of course is
stationary under the action of the flow. The first case is
of particular importance because at certain types of singularity the MCF
is asymptotically self-shrinking.

\bigskip

Before stating our own results, we mention some work that has been
done on the subject: in \cite{AbLa}, all self-shrinking planar
curves
 where classified; in particular, the only simple
self-shrinking curves are the round circles. In \cite{Ang}, the
existence of non spherical self-similar hypersurfaces of
revolution in $\R^n$ were shown; in \cite{An}, we described
rotationally symmetric Lagrangian self-shrinkers and
self-expanders in $\R^{2n}.$ Very recently, a wider class of
self-similar Lagrangian submanifolds has been derived in
\cite{JLT}. On the other hand, few classification results have
been obtained so far. It was shown in \cite{ACR} and \cite{AR}
that the only Lagrangian self-similar submanifolds of $\R^{2n}$
which are foliated by $(n-1)$-dimensional spheres are the examples
found in \cite{An}; in another direction spherical self-shrinkers
have been characterized in \cite{Sm}.

\bigskip

In this note we  give a characterization of the only self-similar
surfaces of $\R^3$ known until now: we first prove that the
self-similar surfaces of revolution discovered by S. Angenent in
\cite{Ang} are the only cyclic self-similar surfaces (Theorem 1),
and next that the cylinders over planar self-similar curves are
the only ruled self-similar surfaces (Theorem 2).

\section{The self-similar equation in coordinates}
Let $X: U \to \R^3$ a local parametrization of some surface
$\Sigma.$ We denote by $E,F$ e $G$ the coefficients of the first
fundamental form of $\Sigma$:
$$ E = |X_s|^2 \quad F = \<X_s, X_t \> \quad G=|X_t|^2.$$
Let $N$ be the unit normal vector given by $N =\frac{X_s \times
X_t}{|X_s \times X_t|}.$ Here and in the remainder of the section,
$\times$ denotes the canonical vectorial product of $\R^3.$ The
coefficients of the second fundamental form are defined to be:
$$ e= \< X_{ss} , N\> \quad f =\< X_{st}, N \> \quad g =\<X_{tt},
N\>.$$
 In order to simplify further calculations, we introduce the
following coefficients, which are proportional to the previous
ones:
$$\bar{e}= \< X_{ss} , X_s \times X_t\>
\quad \bar{f} =\< X_{st}, X_s \times X_t \>
 \quad \bar{g} =\<X_{tt}, X_s \times X_t\>.$$
Rather than the classical formula for the mean curvature,
$$2H = \frac{eG + gE - 2 fF}{EG-F^2} ,$$
 it will be more convenient to use the following one:
$$ 2H= \frac{\bar{e}G + \bar{g}E - 2 \bar{f}F}{(EG-F^2)^{3/2}} .\leqno{(1)}$$

In codimension one, the self-similar equation $\vec{H} + \lambda
X^{\perp}=0$ becomes scalar, namely: $ H + \lambda \<X,N\>=0.$
Moreover, in $\R^3$ have:
$$  \< X, N\> = \frac{1}{\sqrt{EG-F^2}} \< X, X_s
\times X_t \>= \frac{1}{\sqrt{EG-F^2}} \det(X,X_s,X_t).
\leqno{(2)}$$

Finally, from Equations $(1)$ and $(2)$ we deduce:
\begin{lemm}  A surface of $\R^3$ is self-similar if and only if,
for any local parametrization $X: U \to \R^3$ of $\Sigma,$ the
following formula holds:
$$\bar{e}G + \bar{g}E - 2 \bar{f}F + 2 \lambda (EG-F^2) \det(X,X_s,X_t)=0. \leqno{(3)} $$
\end{lemm}

\section{Cyclic surfaces in $\R^3$}

\begin{theo}
Let $\Sigma$ be a self-similar (non minimal) cyclic surface in
$\R^3.$ Then either $\Sigma$ is a round sphere or a surface of
revolution described by S. Angenent (cf \cite{Ang}).

\end{theo}

\begin{lemm}
Let $\Sigma$ be a self-similar cyclic surface in $\R^3.$ Then the
circles of the foliation are parallel or is a piece of a round
sphere.
\end{lemm}

\noindent \textit{Proof of Lemma 2.} The proof is by contradiction
and is based on a method due to J. Nitsche (cf
\cite{Ni1},\cite{Ni2},\cite{Ta}). Let $C(s)$ be a one-parameter
family of circles, $R(s)$ its radius and $\vec{t}(s)$ the unit
normal vector to $C(s).$ There exists some space curve $\ga(s)$
whose unit tangent vector is $\vec{t}(s).$ Moreover, if the
circles are not parallel, the curve $\ga$ is not a straight line,
so its curvature $k(s)$ does not vanish, except in a discrete set
of points. Away from those points, let
$(\vec{t}(s),\vec{n}(s),\vec{b}(s))$ be the Fr\'enet frame related
to $\ga(s).$ Finally, let $z(s)$ be the center of the circle
$C(s).$ Hence, the corresponding cyclic surface is locally
parametrized by
$$ \begin{array}{lccc}  X :
 & I \times \S^1
 &\to&  \R^3\\
&  (s,t)& \mapsto &   R( \vec{n} \cos t + \vec{b} \sin t)+ z.
   \end{array}$$
Following the notation of \cite{Ni1}, we define $(\al, \beta,
\ga)$ to be the coordinates of $z'(s)$ in the Fr\'enet frame
$(\vec{t},\vec{n},\vec{b}).$ A long calculation (cf
\cite{Ni1},\cite{Ta}) shows that $\bar{e}G+\bar{g} E - 2\bar{f} F$
is a trigonometric polynomial in the variable $t$ whose
linearization takes the form:
$$\bar{e}G+\bar{g} E -
2\bar{f} F = \sum_{j=0}^3 a_j \cos(jt) + b_j \sin(jt).$$ For later
convenience, we write the explicit expression of certain of its
coefficients:
$$\begin{array}{ccl}
  a_{3} &=& -\frac{R^{3}k}{2}(k^2R^2+\beta^2-\ga^2), \\
  b_{3} &=& - kR^3\beta\ga, \\
  a_{2} &=& \frac{R^3}{2}(5\al k^2R+\beta'kR-\be k'R-6\beta kR'), \\
  b_{2} &=& \frac{R^3}{3}(\ga'kR-\ga k'R-6\ga kR'). \\
\end{array}$$

\bigskip

We shall now compute the term $(EG -F^2)\det(X,X_s,X_t).$ Firstly
 we define $(p,q,r)$ to be the coordinates of $z(s)$ in the
Fr\'enet frame $(\vec{t},\vec{n},\vec{b}).$ By deriving the
relation $ z= p \vec{t} + q \vec{n} + r \vec{b},$ we get
$$ \left\{ \begin{array}{l} \al = p' - kq, \\ \beta = q' + pk -
\tau r,
  \\ \ga = r' + \tau q.
   \end{array} \right.$$
It follows that the coordinates of $X$ in the frame
$(\vec{t},\vec{n},\vec{b})$ are $(p,q+R \cos t, r + R \sin t).$ We
also have the following expressions for the first derivatives of
the immersion:
$$X_{s}=(\alpha-kR\cos t,\beta+R'\cos t+\tau R\sin t,\gamma-\tau R\cos t+R'\sin t),$$
 $$X_{t}=(0,-R\sin t,R\cos t).$$
Next we calculate:
$$\det(X,X_s,X_t) = \left|  \begin{array}{lll} p & \al - kR \cos t & 0 \\ q + R \cos t &
    \be + R' \cos t + \tau R \sin t  & -R \sin t \\
   r + R \sin t & \ga + R' \sin t - \tau R \cos t & R \cos t \end{array}
   \right|$$
$$ =  \left|  \begin{array}{lll} p & \al - kR \cos t & 0 \\ q + R \cos t &
    \be + R' \cos t  & - R \sin t \\
   r + R \sin t & \ga + R' \sin t  & R \cos t \end{array}
   \right| $$
$$= R \left( p \left| \begin{array}{ll} \be + R' \cos t  & - \sin t \\
    \ga + R' \sin t  & \cos t \end{array}
   \right| + (k R \cos t - \al) \left| \begin{array}{ll} q + R \cos t  & - \sin t \\
    r + R \sin t  &  \cos t \end{array}
   \right| \right)$$
$$ =R \Big( p(\be \cos t + \ga \sin t + R') + (k R \cos t - \al)(q
\cos t + r \sin t + R) \Big)$$
$$ =R \Big( kRq \cos^2 t + kRr \cos t \sin t + (p \be - \al q + kR^2)
\cos t + (p \ga - \al r) \sin t + pR' - \al R \Big)$$
$$ =R \Big( \frac{1}{2} kRq \cos 2 t + \frac{1}{2}kRr \sin 2t + (p \be - \al q + kR^2)
\cos t + (p \ga - \al r) \sin t + \frac{1}{2} kRq+ pR' - \al R
\Big).$$

 The next step is the computation of the first fundamental form:
$$\begin{array}{ccc}
  E &=& \alpha^2+\beta^2+\gamma^2+(R')^2+R^2\tau^2+(2R'\beta -2R\alpha k-2R\gamma\tau) \cos t,\\
    & & +(2R\beta \tau +2R'\gamma) \sin t + R^2 k^2 \cos^2 t ,\\
  F &=& R(-R\tau -\beta \sin t+\gamma \cos t) \\
  G &=& R^2.
\end{array}$$
Thus
$$\frac{EG-F^2}{R^2}=(R^2 k^2 - \ga^2) \cos^2 t - \be^2  \sin^2 t + 2  \be \ga \cos t \sin t $$
$$+ 2(R' \be -  R \al k ) \cos t + 2R'\ga  \sin t
+ (\al^2 + \be^2 + \ga^2 + (R')^2)$$
$$=   \frac{1}{2}(R^2 k^2 - \ga^2 + \be^2 )\cos 2 t + \be \ga \sin
2t +  2(R' \be -  R \al k ) \cos t $$
$$ + 2R'\ga  \sin
t + \left( \frac{1}{2} R^2 k ^2  + \al^2 + \frac{3}{2} \be^2 +
\frac{1}{2} \ga^2 + (R')^2 \right).$$

We deduce that $R^{-3}(EG -F^2)\det(X,X_s,X_t)$ is also a
trigonometric polynomial of order 4, whose linearization take the
form:
$$R^{-3}(EG
-F^2)\det(X,X_s,X_t) =\sum_{j=0}^4 a'_j \cos (jt) + b'_j \sin
(jt).$$
 It follows from Equation 3 that $a'_4$ and $b'_4$ must
vanish, which can be viewed as a linear system in the variables
$q$ and $r$:
$$ \left\{ \begin{array}{cccc}  a'_4=& k \left( \frac{1}{2}(R^2 k^2 - \ga^2 + \be^2)q + \be \ga r  \right) &  = &0 \\
   b'_4= & k \left( - \be \ga q + \frac{1}{2}(R^2 k^2 - \ga^2 + \be^2)r  \right)& =
   &0.
 \end{array} \right. $$
We deduce that either $(i)$ $q$ and $r$ vanish, (which in turn
implies the vanishing of $\ga$), or $(ii)$ the determinant
$\frac{1}{4}(R^2 k^2 - \ga^2 + \be^2)^2 + \be^2 \ga^2$ of the
system vanishes, and then both $R^2 k^2 - \ga^2 + \be^2$ and $\be
\ga$ vanish, so in particular $\be$ or $\ga$ must vanish. But the
vanishing of $\ga$ would imply the vanishing of $R^2 k^2 + \be^2,$
a contradiction (since $R>0$ and $k>0$). So in the second case,
$\be$ vanishes.

\bigskip

\noindent \textbf{First case: $q=r=\ga=0.$}

 We first observe that
here $X(s,t)=R( \vec{n} \cos t + \vec{b} \sin t) + p \vec{t}$ so
that $|X(s,t)|^2 = R^2 + p ^2.$
 Next, using the fact that $\be=pk$ and
$\al=p',$ we get simpler expressions for the following:
$$ \frac{\det(X,X_s,X_t)}{R} =  (p \be  + kR^2)
\cos t + pR' - \al R= k(p^2+R^2) \cos t + pR'-p'R.$$
$$ \frac{EG-F^2}{R^2}=  \frac{1}{2}(R^2 k^2  + \be^2 )\cos 2 t
  + 2(R' \be - R \al k )  \cos t + \left(\frac{1}{2} R^2 k ^2 + \al^2 + \frac{3}{2} \be^2 + (R')^2
\right).$$
$$ =\frac{1}{2} k^2(R^2  + p^2 )\cos 2 t
  + 2k(R'p - Rp' )  \cos t + \left(\frac{1}{2} R^2 k ^2 + (p')^2 + \frac{3}{2} p^2 k^2 + (R')^2
\right).$$
 It follows that:
$$a'_3=\frac{1}{4}(R^2 k^2  + \be^2 )(p \be  + kR^2)=  \frac{1}{4} k^3(R^2+p^2)^2 .$$
On the other hand
 $$ a_{3} = -\frac{R^{3}k}{2}(k^2R^2+ (kp)^2)= -\frac{R^3 k^3}{2}(R^2+p^2).$$
 From Lemma 1, we have
 $$ \sum_{j=0}^4 a_j \cos (jt) + b_j
\sin (jt) +  2\lambda R^3 \left(\sum_{j=0}^4 a'_j \cos (jt) + b'_j
\sin (jt) \right)=0,$$ so that $ \lambda=-\frac{a_3}{2 R^3 a'_3} =
\frac{2}{R^2+p^2}=\frac{2}{|X|^2}.$ It implies that $|X|$ is
constant, thus the surface is a piece of a sphere.


\bigskip
 \noindent
\textbf{Second case:} $\be=0$ and $R^2 k^2 = \ga^2.$

In this case $\ga= \pm Rk$ does not vanish and we have:

$$\frac{\det(X,X_s,X_t)}{R}  = \frac{1}{2} kRq \cos 2 t + \frac{1}{2}kRr \sin 2t + ( - \al q + kR^2)
\cos t + (p \ga - \al r) \sin t + \frac{1}{2} kRq+ pR' - \al R, $$
$$\frac{EG-F^2}{R^2}= - 2 R \al k  \cos t + 2R'\ga \sin t
+ (\al^2 + \ga^2 + (R')^2).$$

As $a_3$ and $b_3$ vanish,
$$ \left\{ \begin{array}{cccc}  a'_3=&  Rk(-qR \al kq -  R' \ga) & = &0 ,\\
   b'_3= & Rk(qR'\ga - r R \al k)  & = &0.
 \end{array} \right. $$
we deduce that either  both $R \al k $  and  $R' \ga$ vanish, or
$q$ and $r$ vanish. We can discard the second case because
$\ga=r'+ \tau q=0,$ a contradiction. Thus, using the fact that
$\al$ and $R' $ vanish, we have
$$\frac{\det(X,X_s,X_t)}{R}  = \frac{1}{2} kRq \cos 2 t + \frac{1}{2}kRr \sin 2t +  kR^2
\cos t + p \ga  \sin t + \frac{1}{2} kRq, $$
$$\frac{EG-F^2}{R^2}=  \ga^2, $$
so that $a_2'= \frac{1}{2}kRq\ga^2$ and
 $b_2'=\frac{1}{2}kRr \ga^2.$
 On the other hand, $a_2$
 and $b_2$ vanish, so by Lemma 1 $a_2'$ and $b_2'$  must vanish as well; again we get a contradiction since
 it implies the vanishing of $q$ and $r.$

\bigskip

\noindent \textit{Proof of Theorem 1.}

By Lemma 2 we know that the circles of a (non spherical)
self-similar cyclic surface must be parallel. Without loss of
generality, we may assume that they are horizontal. Thus the
surface may be locally parametrized by an immersion of the form:
$$ \begin{array}{lccc}  X :
 & I \times \S^1
 &\to&  \R^3\\
&  (s,t)& \mapsto &  (a(s)+ R(s) \cos t,  b(s)+R(s) \sin t, s).
   \end{array}$$
We compute
$$X_s= (a'+ R' \cos t, b' + R' \sin t, 1),$$
$$ X_t = (- R \sin t, R \cos t, 0),$$
from which we deduce the coefficients of the first fundamental
form:
$$ E = (a')^2 + (b')^2 + (R') ^2 + 1 + 2R'( a' \cos t + b'\sin t),$$
$$ F = R(b' \cos t - a'\sin t),$$
$$ G= R^2.$$

We now compute the second derivatives of the immersion:
$$ X_{ss}=(a''+ R'' \cos t, b'' + R'' \sin t , 0),$$
$$ X_{st}= (-R' \sin t, R' \cos t, 0),$$
$$ X_{tt} = (-R \cos t, R \sin t, 0),$$
from which we deduce
$$ \bar{e}= \det(X_{ss} X_s, X_t)= R(- R'' - a'' \cos t + b'' \sin t),$$
$$ \bar{f}=\det(X_{st}, X_s, X_t)=0,$$
$$ \bar{g}= \det(X_{tt}, X_s, X_t)=R^2.$$
Finally, we compute
$$ \det(X,X_s, X_t)= RR's - R^2  + R(a's-a) \cos t + R(b's-b) \sin t.$$
We are now in position to write Equation $(3)$ as a trigonometric
polynomial. There are no terms of order $3$ in $\bar{e}G +
\bar{g}E - 2 \bar{f}F,$ and a straightforward computation shows
that the coefficients in $\cos 3t$ and $\sin 3t $ of $(EG-F^2)
\det(X,X_s,X_t)$  are respectively $(a's-a)((a')^2 +(b')^2)$ and
$(b's-b)((a')^2+(b')^2),$ up to a multiplicative constant. It
follows that either $a'$ and $b'$ vanish, or $a's -a$ and $b's-b$
vanish. The first case is the case of the surfaces of revolution,
which has been treated by S. Angenent in \cite{Ang}. If $a's-a$
and $b's-b$ vanish, we deduce that $a(s)=a_0 s$ and $b(s)=b_0 s,$
for some constants $a_0$ and $b_0.$ It implies the vanishing of
$a''$ and $b''$ and thus the expression  $\bar{e}G + \bar{g}E - 2
\bar{f}F$ becomes a polynomial of degree $1.$ Moreover, the
coefficients in $\cos 2t$ and $\sin 2t $ of $(EG-F^2)
\det(X,X_s,X_t)$  are respectively $(R'-sR)((a_0)^2-(b_0)^2)$ and
$(R'-sR)a_0b_0.$ Again there are two cases: either $R'-sR$
vanishes, or $a_0$ and $b_0$ vanish. If both $a_0$ and $b_0$
vanish, we fall back again on the case of surfaces of revolution.
On the another hand, if $R'-sR$ vanishes, so does
$\det(X,X_s,X_t),$ thus $\bar{e}G + \bar{g}E - 2 \bar{f}F$ must
vanish as well, which means that the immersion is minimal.
Therefore a self-similar cyclic surface must be of revolution and
the proof is complete.

\section{Ruled surfaces in $\R^3$}

\begin{theo}
Let $\Sigma$ be a self-similar ruled surface in $\R^3.$ Then
 $\Sigma$ is a cylinder over a self-similar planar curve.

\end{theo}

\noindent \textit{Proof.} A ruled surface of $\R^3$ may be locally
parametrized by an immersion of the form
$$ \begin{array}{lccc}  X :
 & I \times \R
 &\to&  \R^3\\
&  (s,t)& \mapsto &  \ga(s) t + p(s),
   \end{array}$$
where $|\ga(s)|=1$ and $\<p(s), \ga(s)\>=0.$

 Our discussion being local, we divide the problem in two cases:
either the rulings are parallel, in which case $\ga(s)$ is
constant, or they are not, and then $\ga(s)$ is a regular curve in
$\S^2$ .  The easy task of checking that if the rulings are
parallel, then the ruled surface is a cylinder over a self-similar
planar curve is left to the Reader. Such curves have been
classified by Abresch and Langer (cf \cite{AbLa}). Hence,  we
assume from now on that $\ga(s)$ is a regular spherical curve
 and that $s$ is its arclength parameter. It follows that
 $(e_1, e_2, e_3):=(\ga, \ga', \ga \times \ga')$ is an orthonormal
 frame. Denoting by $k(s)= \<\ga'', \ga \times \ga'\>$ the curvature of $\ga$ in $\S^2,$
  we can write the Fr\'enet equations as follows:
  $e'_1= e_2, e'_2 = ke_3 -e_1 $ and $e_3'= -k e_2.$
Introducing the coordinates of $p$ in the frame $(e_1,e_2,e_3),$
i.e. $p=a e_2 + b e_3,$ we get
$$ p'= -ae_1 + (a'-kb) e_2 + (b'+ ka) e_3.$$

We now compute the first derivatives of the immersion:
$$ X_s= \ga' t + p'  \quad \quad X_t = \ga,$$
from which we deduce the coefficients of the first fundamental
form:
$$ E = t^2  + 2t \< \ga', p'\>+ |p'|^2= t^2 + 2(a'-kb)t + a^2 + (a'-kb)^2 + (b'+ka)^2$$
$$ F = \<\ga, p'\>=a \quad \quad G =1$$
$$ EG-F^2 = t^2 + 2(a'-kb)t +(a'-kb)^2 + (b'+ka)^2.$$

From the second derivatives of the immersion,
$$ X_{ss} = \ga'' t + p'' \quad \quad X_{st} = \ga' \quad \quad
X_{tt}=0,$$ we deduce:
$$ \bar{e}= \det (\ga''t + p'', \ga't +  p', \ga) $$
$$ = t^2 \det (\ga'', \ga', \ga) + t \left(\det(\ga'', p',\ga) + \det(p'', \ga', \ga)
\right) + \det(p'', p', \ga)$$
$$ = k  t^2  + t \left(-k (a'-kb) + \det(p'', \ga', \ga)
\right) + \det(p'', p', \ga),$$
$$ \bar{f} = \det(\ga' , \ga' t + p' , \ga) = \det(\ga', p', \ga)=-(b'+ka),$$
$$ \bar{g} = 0.$$
Finally, we calculate:
$$  \det(X,X_s,X_t)= \det(\ga t + p, \ga't + p', \ga)$$
$$ =   \det( p, \ga' t + p', \ga)$$
$$ =  t \det( p, \ga'  , \ga) + \det(p,p',\ga)$$
$$ = bt + ab'-ba'+ k(a^2+b^2).$$

From Lemma 1 we deduce that the immersion if self-similar if and
only if the following vanishes:
$$ \bar{e}G + \bar{g}E - 2
\bar{f}F  + 2 \lambda (EG-F^2) \det(X,X_s,X_t)=0$$
$$ \Leftrightarrow k t^2  + t \big[ -k(a'-kb) + \det(p'', \ga', \ga)
\big] + \det(p'', p', \ga) - 2 \bar{f}F $$
$$+ 2 \lambda (t^2  + 2(a'-kb)t
+(a'-kb)^2 + (b'+ka)^2 )( bt + ab'-ba'+ k(a^2+b^2))=0.$$

This is a polynomial in $t$ whose coefficient in $t^3$ is $2
\lambda b.$ So $b$ must vanish and we get
$$  k t^2  + t \left[-ka') + \det(p'', \ga', \ga)
\right] + \det(p'', p', \ga) + 2 \lambda \big[t^2  + 2a't +(a')^2
+ (ka)^2 \big] ka^2=0.$$ Now the coefficient in $t^2$ is $k+ 2
\lambda k a^2,$ so either $k$ vanishes, or $\lambda <0$ and $a$ is
a non-vanishing constant.
 If the curvature vanishes, $\ga$ is
a great circle of $\S^2$. As $p = a\ga',$ it follows that the
$X(s,t)=\ga(s) t + a(s) \ga'(s)$ so the image of $X$ lies in the
span of $\ga$ and $\ga'$ and therefore is a piece of a plane. If
$a$ is a constant, using the fact that $p''= -a(1+k^2) e_2 -ak'
e_3,$ we deduce that
$$\det(p'', \ga', \ga) = ak',$$
$$\det(p'', p', \ga)= -a^2 k (1+k^2).$$
Hence the self-similar equation is reduced to
$$  k t^2  + ak' t  - a^2 k (1+k^2) + 2 \lambda (t^2  + (ka)^2 ) a^2 k=0.$$
The coefficient in $t$ is $ak',$ therefore $k$ is constant.
Finally the constant term in the above expression is $a^2 k (1+k^2
+ 2 \lambda k^2 a^2) = a^2 k \big(1+  k^2(1+ 2\lambda a^2)
\big)=a^2 k.$ Again, we get the vanishing of $k,$ hence the
surface is a piece of a plane.


\noindent
Henri Anciaux \\
Departement of Mathematics and Computing \\
Institute of Technology, Tralee \\
Co. Kerry, Ireland \\
henri.anciaux@staff.ittralee.ie \\


\begin{thebibliography}{XXXXX}

 \bibitem[AbLa]{AbLa} U. Abresch, J. Langer, {\em The normalized
 curve shortening flow and homothetic solutions,} J. of Diff.
Geom. \textbf{23} (1986), 175--196

\bibitem[ACR]{ACR} H. Anciaux, I. Castro, P. Romon,
{\em Lagrangian submanifolds of $\R^{2n}$ which are foliated by
spheres,} Acta Math. Sinica (English Series), \textbf{22}(2006)
no. 4, 1197--1214.

\bibitem[An]{An} H. Anciaux, \em Construction of equivariant self-similar solutions
to the mean curvature flow in $\C^n,$ \em Geom. Dedicata,
\textbf{120} (2006), no. 1, 37--48


\bibitem[AR]{AR} H. Anciaux, P. Romon, \em Cyclic and ruled Lagrangian surfaces
 in complex Euclidean space \em math.DG/0703645

\bibitem[Ang]{Ang} S. Angenent, \em Shrinking donuts, \em in
\em Nonlinear diffusion reaction equations \& their equilibrium,
States 3, \em editor N.G. Lloyd, Birkha\"user, Boston, 1992

\bibitem[JLT]{JLT} D. Joyce, Y.-I. Lee, M.-P. Tsui , {\em
  Self-similar solutions and translating solitons for Lagrangian mean curvature flow,}
   arXiv:0801.3721

\bibitem[Ni1]{Ni1} J. Nitsche, \em Lectures on minimal surfaces. Vol. 1
\em Cambridge University Press, Cambridge, 1989

\bibitem[Ni2]{Ni2} J. Nitsche, \em  Cyclic surfaces
 of constant mean curvature, \em
Nachr. Akad. Wiss. Göttingen Math.-Phys. Kl. II 1989,  1--5


\bibitem[Lo]{Lo} R. L\'opez, \em Cyclic surfaces of
constant Gauss curvature,  \em Houston J. Math.  \textbf{27}
(2001), no. 4, 799--805


\bibitem[Sm]{Sm} K. Smoczyk, \em Self-shrinkers of the mean curvature flow in
arbitrary codimension, \em IMRN \textbf{48} (2005), 2983--3004

\bibitem[Ta]{Ta} L. Tavares, Master dissertation, PUC-Rio, 2007



 \end{thebibliography}
\end{document}